\documentstyle{amltd2004}
\begin{document}
\annalsline{158}{2003}
\received{March 13, 2001}
\startingpage{577}
\def\bye{\end{document}}
 \font\tenrm=cmr10
\def\ritem#1{\item[{\rm #1}]}
\def\lra#1{\left\langle#1\right\rangle}
\catcode`\@=11
\font\twelvemsb=msbm10 scaled 1100
\font\tenmsb=msbm10
\font\ninemsb=msbm10 scaled 800
\newfam\msbfam
\textfont\msbfam=\twelvemsb  \scriptfont\msbfam=\ninemsb
  \scriptscriptfont\msbfam=\ninemsb
\def\msb@{\hexnumber@\msbfam}
\def\Bbb{\relax\ifmmode\let\next\Bbb@\else
 \def\next{\errmessage{Use \string\Bbb\space only in math
mode}}\fi\next}
\def\Bbb@#1{{\Bbb@@{#1}}}
\def\Bbb@@#1{\fam\msbfam#1}
\catcode`\@=12

 \catcode`\@=11
\font\twelveeuf=eufm10 scaled 1100
\font\teneuf=eufm10
\font\nineeuf=eufm7 scaled 1100
\newfam\euffam
\textfont\euffam=\twelveeuf  \scriptfont\euffam=\teneuf
  \scriptscriptfont\euffam=\nineeuf
\def\euf@{\hexnumber@\euffam}
\def\frak{\relax\ifmmode\let\next\frak@\else
 \def\next{\errmessage{Use \string\frak\space only in math
mode}}\fi\next}
\def\frak@#1{{\frak@@{#1}}}
\def\frak@@#1{\fam\euffam#1}
\catcode`\@=12
\font\sheaf=eusm10 scaled\magstep1

 \newcommand{\hol}{{{\cal O}}}
\def\Bbb{\bf}
\def\PP{{\Bbb P}}
\def\R{{\Bbb R}}
\def\C{{\Bbb C}}
\def\Z{{\Bbb Z}}
\def\N{{\Bbb N}}
 
\def\eea{\end{eqnarray*}}
\def\bea{\begin{eqnarray*}}

 \title{Moduli spaces\\ of surfaces and real structures}
\shorttitle{Moduli spaces of surfaces and real structures}  
 \acknowledgements{The research of the  author was performed in the realm  of the
  SCHWERPUNKT ``Globale Methode in der komplexen Geometrie",
and of the EAGER EEC Project.}
 \author{Fabrizio Catanese}
\institutions{Mathematisches Institut, Georg-August-Universit\"at,
 G\"ottingen, Germany\\
{\eightpoint {\it Current address\/}:} Universit\"at Bayreuth, Bayreuth, Germany\\
{\eightpoint {\it E-mail address\/}:
Fabrizio.Catanese@uni-bayreuth.de
}}

\centerline{\it This article is dedicated to the memory of Boris  Moisezon}

\vglue26pt

\centerline{\bf Abstract}
\vglue6pt

  We give infinite series  of groups $\Gamma$ and of compact complex
surfaces of general type $S$ with fundamental group $\Gamma$ such that
\begin{itemize}
\item[1)] Any surface $S'$ with
the same Euler number as $S$, and fundamental group $\Gamma$,
is diffeomorphic to $S$.

\item[2)] The moduli space of $S$ consists of exactly two connected
components, exchanged  by complex conjugation.

  Whence,
\begin{itemize}
\item[i)] On the one hand we give simple counterexamples
to the DEF = DIFF question whether deformation type and
diffeomorphism type  coincide for algebraic surfaces.

\item[ii)] On the other hand we get examples of moduli spaces
without real points.

\item[iii)] Another interesting corollary is the existence of complex
surfaces $S$ whose
  fundamental group $\Gamma$
cannot be the fundamental group of a real surface.
\end{itemize}
\end{itemize}

Our surfaces are  surfaces isogenous to a product; i.e.,
  they are quotients
$(C_1 \times C_2)/ G $  of a product of curves by the
free action of a finite group $G$.

They resemble the classical hyperelliptic surfaces, in that
$G$ operates freely on $C_1$, while the second curve is a  {\it triangle
 curve}, meaning that $C_2 / G \equiv \PP ^1$ and the covering
is branched in exactly three points.

\section{Introduction}

Let $S$ be a minimal surface of general type;  then  to $S$
we attach two positive integers $ x = \chi (\hol_S)$, $ y
=K^2_S$ which are invariants of the oriented topological type
of $S$.

The moduli space of the surfaces with invariants $(x,y)$ is a
quasi-projective variety defined over the integers, in
particular it is a real variety (similarly for the Hilbert
scheme of 5-canonical embedded canonical models, of which
the  moduli space is a quotient; cf.\ \cite{Bo}, \cite{Gie}).

For fixed $(x,y)$ we have several possible topological types,
but (by the result of \cite{F}) indeed at most two
if moreover the  surface $S$ is assumed to be simply
connected (actually by [Don1,2], related
results hold  more generally
for the topological types of simply connected compact
oriented differentiable  4-manifolds; cf.\ [Don4,5] for a precise statement, the so-called 11/8 conjecture).

These two cases are distinguished as follows:

\begin{itemize}
\item
  $S$ is EVEN, i.e., its intersection form is even:  
then $S$ is a connected sum of copies of
  $ \PP ^1_{\C} \times \PP ^1_{\C}$ and of a $K3$ surface
if the signature is negative, and of copies of
$ \PP ^1_{\C} \times \PP ^1_{\C}$ and of a $K3$ surface with
reversed orientation if the signature is positive.

\item
$S$ is ODD:   then $S$ is a connected sum of copies of
$\PP ^2_{\C} $ and ${\PP ^2_{\C}}^{\rm opp} .$
\end{itemize}

\numbereddemo{{R}emark}
${\PP ^2_{\C}}^{\rm opp} $ stands for the same manifold as
${\PP ^2_{\C}}$, but with reversed orientation. Beware that
  some authors use the symbol
$\bar{\PP ^2_{\C}} $ for  ${\PP ^2_{\C}}^{\rm opp} $, whereas for us
the notation $\bar{X}$ will denote the conjugate of a complex
manifold $X$ ($\bar{X}$ is just the same differentiable
manifold, but with complex structure $-J$ instead of $J$).
Observe that, if $X$ has odd dimension, then $\bar{X}$
acquires the opposite
orientation of $X$, but if $X$ has even dimension, then
$X$ and $\bar{X}$ are orientedly diffeomorphic.
\enddemo

Recall moreover:

\numbereddemo{{D}efinition}
A real structure $\sigma$ on a complex manifold $X$ is the
{\it datum of an isomorphism} $\sigma:   X \rightarrow \bar{X}$
such that $\sigma^2 = Identity$. One moment's reflection
shows then that $\sigma$ yields an isomorphism between
the pairs $(X, \sigma)$ and $(\bar{X}, \sigma)$.
\enddemo

In general, the fundamental group is a powerful topological
invariant. Invariants of the differentiable structure have been
found by Donaldson, by Seiberg-Witten and several other authors
(cf.\ \cite{Don3}, \cite{D-K}, \cite{Witten}, \cite{F-M3},
\cite{Mor}) and it is well known
that on a connected component of the moduli space the
differentiable structure remains fixed (we use for this result
the slogan DEF $\Rightarrow$ DIFF).

Actually, if two surfaces $S,S'$ are deformation equivalent then
there exists a diffeomorphism carrying the canonical class $K_S$
  $\in H^2 (S, \Z)$
of $S$ to $K_{S'}$; moreover, for minimal surfaces of general type
it was proven  (cf.\ \cite{Witten}, or \cite[Cor.\ 7.4.2, p.\ 123]{Mor})
  that any
diffeomorphism between $S$ and $S'$
  carries $K_S$ either to $K_{S'}$ or to $- K_{S'}$.

Up to recently, the question DEF = DIFF ? was open.

The converse question  DIFF $\Rightarrow$ DEF, asks
whether the existence of an orientation preserving diffeomorphism
between algebraic surfaces $S, S'$ would  imply that $S,S"$ would
be deformation equivalent
(i.e., in the same connected component of the moduli space).
This question was a "speculation" by\break Friedman and Morgan
\cite[p.\ 12]{F-M1} (in the words of the authors, ibidem page~8,
``those questions which we have called speculations$\dots$ seem to require
completely new ideas").

The speculation was inspired by the
successes of gauge theory, and    reading the question I
  thought the answer should be negative, but would not be easy to find.

Recently, (\cite{Man4}) Manetti was able to find
  counterexamples of surfaces with first Betti
number equal to $0$ (but not simply connected).

His result on the one side uses  methods
and results developed in a long sequel of papers
([Cat1,2,3], [Man1,2,3]),
on the other, it uses a rather elaborate construction.

About   the same time of this
paper Kharlamov and Kulikov (\cite{K-K}) gave   a
counterexample for rigid
surfaces, in the spirit of the work of Jost and Yau
([J-Y1,2]).  Here, we have found the following
rather simple series of examples:

\proclaim{Theorem}
Let $S$ be a surface isogenous to a product{\rm ,}  i.e.{\rm ,}  a quotient
$ S = (C_1 \times C_2) /G$ of a product of curves  by the
free action of a finite group~$G$. Then any surface with the
same fundamental group as $S$ and the same Euler number of $S$
is diffeomorphic to $S$. The corresponding moduli space
$ M^{\rm top}_S = M^{\rm diff}_S $ is either irreducible and connected
or it contains
two connected components which are exchanged by complex
conjugation. There are infinitely many examples of the latter
case{\rm ,}  and moreover these moduli spaces are almost all of general type.
\endproclaim

{\it Remark} 1.4.
The last statement is a direct consequence of the results of Harris-Mumford (\cite{H-M}).
\advance\theoremcount by 1

\proclaim{{C}orollary}
{\rm 1) DEF $\neq$ DIFF.}
\vglue4pt
{\rm 2)} There are moduli spaces without real points.
\endproclaim

The more prudent question  of asking whether moduli spaces
with several connected components studied in the previously
cited papers of ours and Manetti would yield diffeomorphic
$4$-manifolds was raised again by Donaldson (in
\cite[pp.\ 65--68]{Don5}), who also illustrated the important
  role played by the symplectic structure of an algebraic surface.
  The referee of this paper points
out an important fact: the standard diffeomorphism between
$S$ and $\overline{S}$ carries the canonical class $K_S$ to
$ - K_{\overline{S}}$, and moreover one could summarize the
philosophy of our topological proof as asserting that there exists
no orientation-preserving self-homeomorphism of $S$, or homotopy
equivalence, carrying $K_S$ to $-K_S$. He then proposes
that one could sharpen the Friedman-Morgan conjecture by asking
whether the existence of a diffeomorphism carrying the canonical class
to the canonical class would suffice to imply deformation
equivalence.

Unfortunately,  this question also  has a negative answer, as
we show in a sequel to this paper (\cite{Cat7}, \cite{C-W}), whose
methods are completely different from the ones of the present paper.

In \cite{Cat7} we give a criterion in order to establish
the symplectomorphism of two algebraic surfaces which are
not deformation equivalent, and show that the examples of
Manetti give a counterexample to the refined conjecture.
Since, however, these examples are not simply connected, we also
discuss some simply connected examples which are not deformation
equivalent: in \cite{C-W} we then show their symplectic equivalence.

  Returning to the examples shown in the present paper,
we deduce moreover, as a byproduct of  our arguments,
  the following:

\proclaim{Theorem}
There are infinite series of groups $\Gamma$ which are fundamental
groups of complex surfaces but which cannot be   fundamental
groups  of a real surface.
\endproclaim

One word about the construction of our examples:
We imitate the hyperelliptic surfaces, in the sense that we
take $ S = (C_1 \times C_2) /G$ where $G$ acts freely on $C_1$,
whereas the quotient $ C_2 /G$ is $\PP ^1_{\C}$.
Moreover, we assume that the projection $\phi: C_2 \rightarrow
\PP ^1_{\C}$ is branched in only three points, namely, we have a
so-called {\it triangle curve}.

It follows that if two surfaces of this sort were
antiholomorphic, then there would be an antiholomorphism of the
second triangle curve (which is rigid).

Now, giving such a branched cover $\phi$ amounts to viewing
the group $G$ as a quotient of the free group with two elements.
Let $a, c$ be the images of the two generators, and set
$ abc = 1$.

We find such a $G$ with the properties that the respective
orders of $a,b,c$ are distinct, whence we show that
an antiholomorphism of
the triangle curve would be a lift of the standard complex
conjugation if the three branch points are chosen to be real, e.g.\
$-1, 0 $ and $+1$.

But such a lifting exists if and only if the group $G$ admits
an automorphism $\tau$ such that $\tau (a) = a^{-1}, \tau (c) =
c^{-1}$.

Appropriate semidirect products do the game for us.

\numbereddemo{{R}emark}
It would be interesting to classify the rigid surfaces,
isogenous to a product, which are not real.
Examples due to Beauville (\cite{Bea}, \cite{Cat6}) yield real surfaces.
\enddemo

\vglue-12pt

\section{A nonreal triangle curve}
\label{first}
\vglue-4pt

Consider  the set $B \subset \PP^1_{\C}$ consisting of  three
real points $ B:   = \{-1, 0, 1\}.$

We choose $2$ as a base point in $\PP^1_{\C} - B$, and   take
the following generators $ \alpha, \beta, \gamma$ of
$ \pi_1 (\PP^1_{\C} - B, 2)$:  
\begin{itemize}
\item
$\alpha$ goes from $2$ to $-1 -\varepsilon$
along the real line, passing through $ + \infty$, then makes
a full turn
counterclockwise around the circumference with centre $-1$ and radius
$\varepsilon$, then goes back to $2$ along the same way on the real line.
\item
$\gamma$ goes from $2$ to $1 +\varepsilon$
along the real line,  then makes a full turn
counterclockwise around the circumference with centre $+1$ and radius
$\varepsilon$, then goes back to $2$ along the same way on the real line.
\item
$\beta$ goes from $2$ to $1 +\varepsilon$
along the real line,   makes a half turn
counterclockwise around the circumference with centre $+1$ and radius~$\varepsilon$, reaching $1 -\varepsilon$,
then proceeds along the real line  reaching $+ \varepsilon$,
  makes a full turn
counterclockwise around the circumference with centre $0$ and radius
$\varepsilon$,  goes back to $1 -\varepsilon$ along the same way
on the real line, makes again  a half turn
clockwise around the circumference with centre $+1$ and radius~$\varepsilon$, reaching $1 +\varepsilon$; finally
it proceeds along the real line returning to~$2$.
\end{itemize}

An easy picture shows that $ \alpha,  \gamma$ are  free generators
of $ \pi_1 (\PP^1_{\C} - B, 2)$ and
$$ \alpha \beta \gamma = 1 .$$

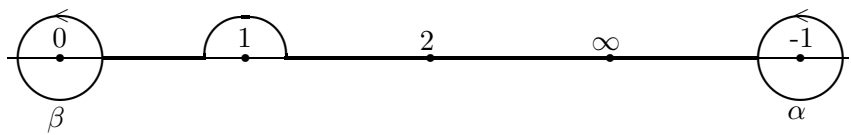
\font\mate=cmmi8
\def\min{\hbox{\mate \char60}\hskip1pt}
\begin{center}
\begin{picture}(300,55)(0,-30)
%
%
\put(0,0){\line(1,0){300}}
\put(238,0){\line(1,0){80}}
%
%
\put(20,0){\circle*{3}}
\put(17,4){0}
\put(90,0){\circle*{3}}
\put(87,4){1}
\put(160,0){\circle*{3}}
\put(156,3){2}
\put(228,0){\circle*{3}}
\put(221,3){$\infty$}
\put(300,0){\circle*{3}}
\put(296,4){-1}
{\thicklines
%
%
\put(20,0){\circle{30}}
\put(36,0){\line(1,0){39}}
\put(90,0){\oval(31,31)[t]}
\put(105,0){\line(1,0){55}}
\put(17,14){$\min$}
\put(15,-26){$\beta$}
%
%
\put(300,0){\circle{30}}
\put(160,0){\line(1,0){124}}
\put(297,14){$\min$}
\put(295,-23){$\alpha$}
}
\end{picture}
\end{center}

With this choice of basis,  we have provided an isomorphism
of\break $ \pi_1 (\PP^1_{\C} - B, 2)$ with the group
$$  T_{\infty}:   =\langle \alpha, \beta, \gamma | \
  \alpha \beta \gamma = 1 \rangle .$$

For each finite group $G$ generated by two elements $a,b$,
passing from Greek to italic letters we obtain a tautological
surjection
$$ \pi:   T_{\infty}  \rightarrow G .$$
That is, we set $\pi (\alpha) = a, \pi (\beta) = b$ and we define
$ \pi (\gamma):  = c.$ (then $abc=1$).

\numbereddemo{{D}efinition}
We let the triangle curve $ C$ associated to $\pi$
be the {\it Galois covering} $f:  C \rightarrow
\PP^1_{\C}$,  branched on $B$  and with group $G$  determined by the
chosen isomorphism $ \pi_1 (\PP^1_{\C} - B, 2) \cong T_{\infty}$
and by the group epimorphism $\pi$.
\enddemo

\numbereddemo{{R}emark}
Under the above notation, we set $m,n,p$ the periods of the
respective elements $a,b,c$ of $G$ (these are the branching
multiplicities of the covering $f$). Composing $f$ with a projectivity
we can assume that $m \leq n \leq p$.
\enddemo

Notice that the Fermat curve $ C:  = \{(x_0,x_1,x_2) \in \PP^2_{\C} |
x_0^n + x_1^n + x_2^n = 0\}$ is in two ways a triangle curve,
since we can take the quotient of $C$ by the group $G:  = (\Z/n) ^2$ of
diagonal projectivities with entries n-th roots of unity, but
also by the full group $A= {\rm Aut}(C)$ of automorphisms, which
is a semidirect product of the normal subgroup $G$ by the symmetric
group exchanging the three coordinates.
For $G$ the three branching multiplicities are all equal to $n$,
whereas for $A$ they are equal to $(2,3,2n)$.

Another interesting example is provided by the {\it Accola curve}
(cf.\ [ACC1,2]),
the curve $ Y_g $ birational to the affine curve of equation
$$ y^2 = x^{2g+2} - 1 .$$

If we take the group $G \cong \Z/ 2 \times \Z/(2g+2)$ which acts
multiplying $y$ by $-1$, respectively $x$ by a primitive $2g+2$-root
of $1$, we realize $ Y_g $ as a triangle curve with
  branching multiplicities $ (2, 2g+2, 2g+2)$. However, $G$ is not 
the full automorphism group; in fact if we add the transformation
sending $x$ to $1/x$ and $y$ to $iy/x^{g+1}$, then we get a nonsplit extension of $G$ by $\Z/2$ (which is indeed
the full group of automorphisms of $ Y_g $ as   is well known and as also
follows from the
next lemma), a group which represents $ Y_g $ as a triangle curve with
  branching multiplicities $(2, 4,  2g+2)$.

One can get many more examples by taking unramified coverings
of the above curves (associated to characteristic subgroups
of the fundamental group).

The following natural question arises then:
which are the curves which admit more than one realization as
triangle curves?

We are not aware whether the answer is already known in the literature,
but (although this is not strictly needed for our purposes)
we will show in the next lemma  that this situation
is rather exceptional
if the branching multiplicities are all distinct:

\proclaim{Lemma}
Let $f:   C \rightarrow \PP^1_{\C} = C/G$ be a triangle covering
where the branching multiplicities $m,n,p$ are all distinct
{\rm (}\/with the assumption that $m < n < p$\/{\rm).}
The group $G$ equals the full group $A$
of automorphisms of $C$ if the triple is not $(3, m_1, 3 m_1)  $
  or $(2, m_1, 2 m_1)  $.
\endproclaim

\demo{Proof} I. By Hurwitz's formula the cardinality of $G$ is in general
given by the formula $$ |G| = 2 (g-1) (1-1/m-1/n-1/p)^{-1} .$$
 
II. Assume that $A\neq G$ and let $F: \PP^1_{\C} = C/G
\rightarrow \PP^1_{\C} = C/A$ be the induced map.
Then $f': C \rightarrow \PP^1_{\C} = C/A$ is again a triangle covering,
otherwise the number of branch points would be $\geq 4$ and
we would have a nontrivial family of such Galois covers
with group $A$ (the cross ratios of the branch points
would provide locally nonconstant holomorphic functions on the
corresponding subspace of the moduli space).
Whence, also a nontrivial family of $G$ -covers,
a contradiction.
\vglue4pt
III. Observe that, given two points $y,z$ of $C$, $f'(y) = f'(z) $
if and only if $z \in  A y $ and then the branching indices  of
$y,z$ for $f'$ are the same. On the other hand, the branching
index of $y$ for $f'$ is the product of the branching
index of $y$ for $f$ times the one of $f(y)$ for $F$.
\vglue4pt
IV. We claim now that the three branch points of $f$
cannot have distinct images
through $F$: otherwise the branching multiplicities $m' \leq n'
\leq p'$ for $f'$ would be not less than the respective multiplicities
for $f$, and by the analogue of formula I for $|A|$ we would
obtain $ |A| \leq |G|$, a contradiction.
\vglue4pt
V. Note that if the branching multiplicities $m,n,p$ are all
distinct, then $G$ is equal to its normalizer in $A$,
because if $\phi \in A$, $ G  = \phi G \phi ^{-1}$,
  then $\phi$ induces an automorphism
of $\PP^1_{\C}$, fixing $B$, and moreover such that it sends
each branch point to a branch point  of the same order.
Since the three orders are distinct, this automorphism
must be the identity on $\PP^1_{\C}$, whence $\phi \in G.$
 \vglue4pt
VI. Let $x_1, x_2, x_3$ be the  branch points of $f$ of respective
multiplicities $m_1,m_2,m_3$ (that is, we consider again the
three integers $m,n, p$, but allow  another ordering).
Suppose now that $F(x_1)=F(x_2) \neq F(x_3)$:
  we may clearly assume $m_1 < m_2$.
  Thus the branching multiplicities for $f'$ are $n_0, n_2,
n_3$, where $n_2,  n_3$  are the respective multiplicities of
$F(x_2) \neq F(x_3)$.
Thus $n_2$ is a common multiple of $m_1,m_2$, $n_2 = \nu_1 m_1 =
\nu_2 m_2$,  $n_0$ is greater or equal to $2$,
  $n_3 = m_3 \nu_3 $,
whence $m_2 \leq n_2,  n_2 \geq 2 m_1$.

We obtain 
\begin{eqnarray*}
|A|/|G|& \leq&  (1- 1/m_3- 3/n_2)
(1- 1/2 -1/m_3- 1/n_2)^{-1}  \\
&=&  2 \frac{ n_2 m_3 - n_2 - 3m_3 }
{n_2 m_3 - 2 n_2 - 2 m_3 }= 2 +\frac{ 2 n_2 - 2 m_3 }
{n_2 (m_3 - 2)- 2 m_3 }.
\end{eqnarray*}
Thus $ |A|/|G| \leq  2$ if $m_3 \geq 5$,
$ |A|/|G| \leq  3$ if $m_3 = 4$.
 \vglue4pt
VII. However, if $ |A|/|G| \leq  2$ then $G$ is normal
in $A$; thus, by our assumption and by V, $G=A$.
Now, we need only to take care of the possibility
$ |A|/|G| \geq  3$.
 \vglue4pt
VIII. Under the hypothesis of VI, we get $ d: = {\rm deg} (F) =
  |A|/|G| = k_0 n_0 $. Since $n_0 \geq 2$, if $d=3$ we get
$n_0 = 3$. We also have
$$ {\rm (i)}\quad  d =  \nu_3 (1 + k_3 m_3),\qquad
{\rm (ii)}\quad  d = \nu_1 + \nu_2 +  k_2 n_2   \  \ (k_2, k_3 \geq 0).$$

  Now, if $m_3 = 4$ we get $d= 3 = n_0 = \nu_3$; but then $F$
cannot have
further ramification points, contradicting $\nu_1 \geq 2$.

If instead $m_3 = 3$ the above inequality yields
$d = |A|/|G| \leq 3 + n_2/ (n_2 - 6)$. But $n_2 = \nu_1 m_1 \geq 8$
(this is obvious if  $m_1 \geq 4$, else $m_1 =2$ but then $m_2 \geq
8$).

Next, $n_2 \geq 8$ implies $ d \leq 7$. From  (ii)  and  $n_2 \geq 8$
follows then  $k_2 = 0$, whence  $d = \nu_1 + \nu_2 $.

Then the previous inequality yields
$$  d \leq 2 \frac{2 n_2 - 3d}{n_2 -6 };  $$
i.e., $d(n_2 -6) \leq 4 n_2 - 6d  $,  whence $d \leq 4$.

If $d=3$ we get the same contradiction from $d = n_0 = \nu_3$.
Else, $d=4$ and equality holds, whence $\nu_3 = 1$,$n_0=2$,
and  $\nu_1 = 3, \nu_2 = 1$.
In this case we get $d = |A|/|G| = 4$, $m_3 = n_3 = 3 $, $n_0 = 2$,
$ n_2 = 3 m_1 = m_2 \geq 8$.
Then the branching indices are
$$ (3, m_1, 3 m_1)  \ {\rm for \ G  \ and } \ (2, 3,  3 m_1) {\rm \ 
for \ A}.$$

Assume finally that $m_3 =  2$. If $n_3=2$, then $n_0 \geq 3$,
thus the usual inequality gives
$$ d \leq  (1/2 - \frac{\nu_1 + \nu_2}{n_2}) \frac{6 n_2}{n_2 -6} =
3 \frac{n_2 - 2(\nu_1 + \nu_2)}{n_2 -6} \leq 3.$$

But again $d = 3$ implies $n_0 = 3$, and $\nu_3 = 3$ yields
the usual contradiction.
Thus  $\nu_3 = 1 =  \nu_2 $ and then $ m_3 = n_3 = 2$, $\nu_1 = 2$, $n_0 = 3$,
$ n_2 = 2 m_1 = m_2 $ and we have therefore the case $d =3$ and
branching indices
$$ (2, m_1, 2 m_1)  \ {\rm for \ G  \ and } \ (2, 3,  2 m_1) {\rm \ 
for \ A}.$$
\vglue4pt
IX. There remains the case where $F(x_1)=F(x_2) =
F(x_3)$. Then the branching order of $f'$ in $F(x_i)$
is a common multiple $\nu$ of $m,n,p$, and we get the
estimate
\begin{eqnarray*}
 |A|/|G| &\leq & (1- 1/m- 1/n -1/p)
(1- 1/2 -1/3- 1/\nu)^{-1} \\
&= &(1- 1/m- 1/n -1/p)
\frac{6\nu}{\nu - 6} . 
\end{eqnarray*}

Now, if $p < \nu$, then $ \nu \geq 2p, \nu \geq 3n,\nu \geq 4m$;
thus $ |A|/|G| \leq
\frac{6 (\nu-9)}{\nu - 6} < 6. $
However, looking at the inverse image of $F(x_i)$ under $F$,
we obtain $$(*) |A|/|G| \geq  \nu/m + \nu/n +\nu/p,  $$
whence $ |A|/|G| \geq 9$, a contradiction.

Thus $p=\nu$, and then from this equality follow also the
further inequalities $  \nu \geq 2n,\nu \geq 3m$.
We get  $ |A|/|G| \leq  6 $ from the first inequality,
and from $(*)$ we derive that $ |A|/|G| \geq  6 $.

The only possibility is:   $ |A|/|G| = 6 $, $p=3m, p=2n$.

In this case therefore the three local monodromies of $F$
are given by
three permutations in six elements, with  cycle decompositions
of respective types $(1,2,3), (n)^{k},(n')^{k'}$, where
$ nk = n'k' = 6$.
The Hurwitz formula for $F$\break (${\rm deg}\ F =6$) shows that the respective
types must then be $(1,2,3), (2,2,2), (3,3)$.

We will conclude then, deriving a contradiction
by virtue of the following Lemma.

\proclaim{Lemma}
Let $ \tau, \sigma$ be permutations in six elements of respective
types $(2,2,2), (3,3)$. If their product $\sigma \tau$
has a fixed point{\rm ,} then it has a cycle decomposition
of type $(1,4,1)$.
\endproclaim

\demo{Proof}
We will prove the lemma by suitably labelling the six elements.
Assume that $2$ is the element fixed by $\sigma \tau$: then we
label $1: = \tau (2)$.
Since   $\sigma (1) =2$, we also label $ 3: = \sigma (2)$.
Further we label  $4: = \tau (3)$, $5: = \sigma (4)$,
so that $\tau$ is a product of the three transpositions
$(1,2)$, $(3,4)$, $(5,6)$, while $\sigma$ is the product of
the two three-cycles $(1,2,3)$, $(4,5,6)$.

An easy calculation shows that $\sigma \tau$ is the four-cycle
$(1,3,5,4)  $. 
\enddemo

  \numbereddemo{{R}emark}
The above proof of lemma 2.3 provides explicitly a realization of $ T:   =
T(3, m_1, 3 m_1)  $ as a (nonnormal) index $4$ subgroup of
$ T':   =  T(2,3, 3 m_1)  $, resp.\ of $ T:   =
T(2, m_1, 2 m_1)  $ as a (nonnormal) index $3$ subgroup of
$ T':   =  T(2,3, 2 m_1)  $.

For every finite index normal subgroup $K$ of $T'$, with $ K \subset T$, we
get $   G:  = (T / K) \subset  A:  = (T' /K)  $ and corresponding 
triangle curves.

Thus the exceptions can be characterized.
\enddemo

We come now to our particular triangle curves.
Let $r, \ m$ be positive integers  $r\geq 3,  m \geq 4$ and set
$$ p:  = r^m -1,\  n:= (r-1) m \ . $$
Notice that the three integers $m < n < p $ are distinct.

Let $G$ be the following semidirect product of $ \Z/p $ by $ \Z/m$:
$$  G:   =\langle  a,c \ |\  a^m = 1,  \ c^p = 1, \ a c a^{-1} = c^r \rangle $$

The definition is well posed (i.e., the semidirect product of
$ \Z/p $ by $ \Z$ given by
$G'':   = \lra{ a,c \ |\ c^p = 1, \ a c a^{-1} = c^r } $ descends to
a semidirect product of $ \Z/p $ by $ \Z/m$) since
$$ a^i c a^{-i} \ = c^{r^i} $$
and, by very definition of $m$, $r^m \equiv 1 \ ({\rm mod} \ p) .$

\proclaim{Lemma}
Define $ b \in G$ by  the equation $ abc = 1$. Then the period of
$b$ is exactly $n$.
\endproclaim

\demo{Proof}
The elements of $G$ can be uniquely written as
$\{ c^i a^j | \ 0 \leq i \leq p-1,  \ 0 \leq j \leq m-1 \}.$
The period of $b$ equals the one of its inverse, namely, $ ca$.
Now,
$$    (ca)^i =  c \ (aca^{-1}) \ (a^2ca^{-2}) \cdots
(a^{i-1}ca^{-(i-1)}) \ a^i \ = \\
  c ^{ 1 + r + .. r^{i-1}} \ a^i =  c^{ \frac{r^i - 1 }{r-1}}
\ a^i.$$
Whence, $b^i = 1 $ if and only if $ m | i$ and
$ p |\frac{r^i - 1 }{r-1}$.

Let therefore $ mk $ be the period of $b$; then $k$ is the smallest
integer with $\frac{r^{mk} - 1 }{r-1} \equiv 0 \ {\rm mod} \ (r^m -1).$
Since $\frac{r^{mk} - 1 }{r-1} = \frac{r^{mk} - 1 }{r^m -1}
\frac{r^{m} - 1 }{r-1}$,  all we want is $ \frac{r^{mk} - 1 }{r^m -1}
\equiv 0 \ ({\rm mod} (r-1))$; however, $ \frac{r^{mk} - 1 }{r^m -1}
\equiv  \sum_{j=0}^{k-1} r^{mj}  \equiv k \ ({\rm {\rm mod }\ }(r-1)).$
Therefore $ k = r-1$ and the period of $b$ equals $n$.
\enddemo

\proclaim{Proposition}
The triangle curve $ C$ associated to $\pi$ is not antiholomorphically
  equivalent to itself {\rm (}\/i.e.{\rm ,} it is not isomorphic to its conjugate\/{\rm).}
\endproclaim

\demo{{P}roof}
We shall derive a contradiction assuming the existence of
   an antiholomorphic automorphism $\sigma$ of $C$.

\demo{Step I} $G$  = $A$, where $A$ is the group
of holomorphic automorphisms of~$C$, $A: = {\rm Bihol}(C,C)$.
\enddemo

{\it Proof}. This follows from the previous Lemma 2.3, since in this case
we assumed $m \geq 4$, and since $p=r^m-1,  n =m(r-1)$,
obviously $$p  >  (r-1 +m)(r-1)^{m-1} \geq (2+m) 9 (r-1) >  2n.$$

\demo{Step II} If $\sigma$ exists, it must be a lift of complex
conjugation.
\enddemo

{\it Proof}. In fact $\sigma$ normalizes ${\rm Aut}(C)$, whence it must induce
an antiholomorphism of $\PP^1_{\C}$ which is the identity on
$B$, and therefore must be complex conjugation.

\demo{Step III} Complex conjugation does not lift.
\enddemo 

{\it Proof}. This is purely an argument about covering spaces:
complex conjugation acts on $ \pi_1 (\PP^1_{\C} - B, 2) \cong
T_{\infty} $,
as   is immediate   with our choice of basis,
by  the automorphism $\tau$ sending $\alpha$,  $\gamma$
to their respective inverses.

Thus, complex conjugation lifts if and only if $\tau$ preserves
the normal subgroup $K:   = {\rm ker} (\pi)$. In turn, this means that there
is an automorphism $ \rho:   G \rightarrow G$ with
$$ \rho (a) = a^{-1}, \  \rho (c) = c^{-1}. $$

Recall now the relation $ a c a ^{-1} = c^r $:  
Applying $\rho $,
we would get $ a^{-1} c ^{-1}a\break  = c^{-r }$, or, equivalently,
$$\pagebreak a^{-1} c a  = c^{r }.$$

But then we would get $ c = a ^{-1} (a c a ^{-1}) a = c^{r^2}  $, 
which holds only if
$$ r^2 \equiv 1 \ ({\rm mod \ }p).$$
This is the desired  contradiction,
because $ r^2 - 1 < r^m - 1 = p $.
\hfill\qed\enddemo
\vfill
{\it Remark} 2.8.
What the above proposition says can be rephrased in the following terms
(cf.\ \cite[pp.\ 29--31]{Cat6}).

Denote by $\Pi_g$ the fundamental group of a compact curve of
genus $g$.

The epimorphism $\pi:   T_{\infty}  \rightarrow G $  factors through
an epimorphism $\pi'$ of the triangle group
$ T(m,n,p):   =  \lra{a,b,c \ |\  a^m = 1, \  b^n =1,
\ c^p = 1, \ a b c  = 1}$ onto $G$, and once an isomorphism
$ {\rm ker}\ \pi' \cong \Pi_g $ is fixed, the pair
$(C, \ {\rm ker}\ \pi' \cong \Pi_g)$ yields a point
in the Teichm\"uller space ${\cal T}_g$. This point is the
only fixpoint for the action of $G$ on ${\cal T}_g$ induced
by the natural homomorphism $ G:   \rightarrow {\rm Out}
({\rm ker}\ \pi' \cong \Pi_g)$.
The pair $(\bar{C},{\rm ker}\ \pi' \circ \tau') $ corresponds to the epimorphism
$\pi' \circ \tau':   T_{m,n,p}  \rightarrow G .$
What we have shown is that $\bar{C}, C$ yield different points
in the moduli space ${\cal M}_g = {\cal T}_g / {\rm Out} \ (\Pi_g).$
Thus, $\bar{C}$ and $C$ correspond to two topological actions of $G$ which are conjugate
by an orientation reversing homeomorphism, but not by an orientation-preserving one.
\enddemo
 
\vglue-12pt
\section{Theorems, and corrigenda}
\label{second}
\vglue-4pt

We begin this section by recalling some results of (\cite{Cat6}),
and we draw some consequences for real surfaces.
For  one of the theorems of (\cite{Cat6}) we shall need to make
  a small correction which,
although it amounts only  to  remembering that $ (-1)^2 = 1  $,
will be completely crucial to our argument.

Recall  \cite[3.1--3.13]{Cat6}:

\numbereddemo{{D}efinition}
A {\it projective surface} $S$ is said to be {\it isogenous} to a  (higher) product
if it admits a finite unramified covering by
  a product of curves of genus $\geq 2$. In this case, there
exist  Galois realizations $ S = (C_1 \times C_2)/ G$, and  each
such Galois realization dominates a uniquely determined minimal
one. Note that $S$ is said to be of nonmixed type if $G$ acts via a product
action of two respective actions on $C_1, C_2$. Otherwise $S$ is
of mixed type and it has a canonical unramified double cover
which is of unmixed type  and $C_1 \cong C_2$ (see 3.16 of
\cite{Cat6} for more details on the realization of surfaces of mixed
type). In the latter case the canonical double cover corresponds
to a subgroup $ G^0 \subset G$ of index $2$.
\enddemo

\proclaim{Proposition}
Let $S, S'$ be surfaces isogenous to a higher product{\rm ,} and let
$\sigma:   S \rightarrow S'$ be an antiholomorphic isomorphism.
Let moreover $ S = (C_1 \times C_2)/ G${\rm ,}
$ S' = (C'_1 \times C'_2)/ G'$ be the respective
minimal Galois realizations. Then{\rm ,} up to a  possible exchange of
$C'_1 $ with $ C'_2 ${\rm ,} there exist antiholomorphic isomorphisms
$\tilde{\sigma}_i,  \ i= 1,2 $ such that $\tilde{\sigma}:  =
\tilde{\sigma}_1 \times \tilde{\sigma}_2  $ normalizes the action
of $G${\rm ;} in particular $\tilde{\sigma}_i $ normalizes the action
of $ G^0 $ on $C_i$. \pagebreak
\endproclaim

\demo{Proof}
Let us view $\sigma$ as yielding a complex isomorphism
$\sigma:   S \rightarrow \bar{S'}$.
Consider the exact sequence corresponding  to the minimal
Galois realization
$ S = (C_1 \times C_2)/ G$,
$$1 \rightarrow H:   = \Pi_{g_1} \times \Pi_{g_2}
\rightarrow  \pi_1(S) \rightarrow G \rightarrow 1. $$

Applying $\sigma_*$  to it, we infer by Theorem 3.4 of (\cite{Cat6})
that there is an exact sequence associated to a Galois
realization of $ \bar{S'}$. Since $\sigma$ is an isomorphism,
we get a minimal one, which is however unique.
Whence, we get an isomorphism    $ \tilde{\sigma}:  
(C_1 \times C_2)  \rightarrow (\bar{C'_1} \times \bar{C'_2}) $,
which is of product type by the rigidity lemma (e.g., Lemma 3.8
of \cite{Cat6}). Moreover, this isomorphism must normalize
the action of $ G \cong G'$, which is exactly what we claim.
\enddemo

The following is the correction of Theorems 4.13,  4.14
of \cite{Cat6}:

\proclaim{Theorem}
Let $S$ be a surface isogenous to a product{\rm ,} i.e.{\rm ,} a quotient
$ S = (C_1 \times C_2) /G$ of a product of curves  by the
free action of a finite group~$G$. Then any surface $S'$ with the
same fundamental group as $S$ and the same Euler number of $S$
is diffeomorphic to $S$. The corresponding moduli space
$ M^{\rm top}_S = M^{\rm diff}_S $ is either irreducible and connected
or it contains
two connected components which are exchanged by complex
conjugation.
\endproclaim

\demo{Proof}
The only modification in the proof given in \cite{Cat6} occurs
on the last lines of page 30.

As in the previous proposition, an isomorphism
between the fundamental groups of $S$, resp.\ $S'$ yields
  a differentiable action of $G$ on the product of curves
$(C'_1 \times C'_2) $ yielding  the minimal
Galois realization of $S'$. In fact the above isomorphism
of fundamental groups,  by unicity of the Galois realization,
yields an isomorphism of $H$ with $H'$. This isomorphism
yields an orientation-preserving diffeomorphism
$ (C_1 \times C_2) \rightarrow
(C'_1 \times C'_2) $ which is of product type.

Now, the diffeomorphisms between the
respective factors are either both orientation-preserving
(this was the case we were considering in the argument in loc.\ cit.),
or both orientation-reversing.

In the latter case,   the topological action
of $G^0$ on the product of the conjugate curves
$(\bar{C'_1} \times \bar{C'_2}) $, which is of product type,
yields actions  of $G^0$  on the respective factors
$\bar{C'_1},  \bar{C'_2}  $ which are
  of the same oriented topological type as the respective actions on $ 
C_1,  C_2 $
(again here we might have to exchange the roles of
$C'_1,  C'_2  $ if the genera $g_1, g_2$ are equal).
Therefore we conclude in this case that the conjugate
of the surface $S'$ belongs to the irreducible subset of
the moduli space containing $S$.
\enddemo

We are now going to explain  the construction
of our examples:  

Let $G$ be the semidirect product group   constructed in Section~$2$, and let $C_2$ be the corresponding triangle
curve.

By the formula of Riemann Hurwitz the genus of $C_2$ equals
$$ g_2 =  1 + \frac{1}{2} [(m-1)(r^m-2) -1 - \frac{r^m-1}{r-1}].$$

Let moreover $g'_1$ be any number greater than or equal to $2$,
and consider the canonical epimorphism $\psi$ of $\Pi_{g'_1}$
onto a free group of rank $g'_1$, such that in terms of the standard
bases $ \ a_1,  b_1, .... a_{g'_1},\ b_{g'_1} $, respectively
$ \gamma_1,.. \gamma_{g'_1}$,
we have
$\psi (a_i) = \psi(b_i) = \gamma_i$.
Then compose   $\psi$  with any epimorphism of the free group
onto $G$, e.g. it suffices to compose with any $\mu$ such that
$ \mu (\gamma_1) = a$,  $\mu (\gamma_2) = b$ (and
$\mu (\gamma_j)$ can be chosen to be  whatever we want for $j \geq 3$).

For any point $C'_1$ in the Teichm\"uller space we obtain a canonical
covering associated to the kernel of the epimorphism $ \mu \circ
\psi:   \Pi_{g'_1} \rightarrow G$; call it $C_1$.

\numbereddemo{{D}efinition}
Let $S$ be the surface $S:   =(C_1 \times C_2) / G$ ($S$ is smooth
because $G$ acts freely on the first factor).
\enddemo

\proclaim{Theorem}
For any two choices  $C'_1(I),  C'_1(II) $
of $C'_1$ in the\break Teichm{\rm \"{\it u}}ller space there are surfaces  $S(I), S(II)$
such that $S(I)$ is never isomorphic to $\bar{S}(II)$.
When $C'_1$ is varied{\rm ,} there is a connected component of the moduli space{\rm ,}
which has only one other connected component{\rm ,} given by the conjugate of
the previous one.
\endproclaim

\demo{Proof}
By Theorem 3.3 it suffices to show the first statement, because
we know already, by the rigidity of the second triangle curve,
that  we get a connected component of the moduli space varying $C'_1$.
By Proposition 3.2, it follows the fact that if $S(I)$ were isomorphic to $\bar{S}(II)$,
then there would be an antiholomorphism of $C_2$ to itself.
This is however excluded by Proposition~2.3.
\enddemo

We come now to the last result:

\proclaim{Theorem}
Let $S$ be a surface in one of the families constructed above.
Assume moreover that $X$ is another complex surface such
that $\pi_1(X) \cong \pi_1(S)$. Then $X$ does not admit any real
structure.
\endproclaim

\demo{Proof}
Observe that since $S$ is a classifying space for the fundamental
group of $\pi_1(S)$, then by the isotropic subspace theorem of
 \cite{Cat3} 
  the Albanese mapping of $X$ maps onto a curve $C'(I)_2$
of the same genus as $C'_2$.

Consider now the unramified covering $\tilde{X}$ associated to
the kernel of the epimorphism $\pi_1(X) \cong \pi_1(S) \rightarrow G$.

Again by the isotropic subspace theorem, there exists a holomorphic
map $\tilde{X} \rightarrow C(I)_1 \times C(I)_2$, where moreover
the action of $G$ on $\tilde{X}$ induces  actions of $G$ on
both factors which either have the same oriented topological types as 
the actions
of $G$ on $C_1$, resp.\ $C_2$,  or have both the oriented topological types
of the actions on the respective conjugate curves.

By the rigidity of the triangle curve $C_2$, in the former case
$ C(I)_2 \cong  C_2$, in the latter $ C(I)_2 \cong  \bar{C_2}$.

Assume  now that $X$ has a real structure $\sigma$: then the same
argument as
in [C-F, \S 2] shows that $\sigma$ induces  a product
antiholomorphic map $\tilde{\sigma}:   C(I)_1 \times C(I)_2
\rightarrow C(I)_1 \times C(I)_2$. In particular, we get a
nonconstant antiholomorphic map of $C_2$ to itself,
contradicting Proposition 2.3 .
\enddemo

Finally, in  \cite{Cat6}, motivated by  some examples by Beauville
(\cite[p.\ 159]{Bea}) we gave the following:

\numbereddemo{{D}efinition}
A {\it Beauville surface} is a rigid surface isogenous to a\break product.
\enddemo

However, in commenting in five lines on where  the problem of
the classification of such surfaces lies, I confused together the
nonmixed type and the mixed type (which is more difficult to get).

Therefore,  I would simply like here to comment that
  to obtain a Beauville
surface of nonmixed type it is equivalent to give a finite group
$G$ together with two systems of generators $\{a,b\}$ and
$\{a',b'\}$ which satisfy a further property, denoted
by  (*) in the sequel.

In fact, the choice of the two systems of generators yields two
epimorphisms  $\pi, \pi':   T_{\infty}  \rightarrow G $ where we
recall that
$  T_{\infty}:   = \lra{\alpha, \beta, \gamma | \
  \alpha \beta \gamma = 1 }$
is the fundamental group of the projective line minus three points.

We get corresponding curves $C, C'$ with an action of $G$, and
the product action of $G$ on $ C \times C'$ is free if and only
if,  defining $c,c'$ by the properties $abc= a'b'c' = 1$,
and letting $\Sigma$   be the union of the conjugates of the cyclic
subgroups generated by $a,b,c$ respectively, and
defining $\Sigma'$
analogously,   we have
$$ (*) \  \Sigma \cap \Sigma' = \{ 1_G \} .  $$

In the mixed case, one requires instead that the two systems
of generators be related by an automorphism $\phi$ of $G$
which should satisfy the further conditions:

\begin{itemize}
\item
i) $\phi^2$ is an inner automorphism, i.e., $\phi^2 = {\rm Int}_{\tau}$
for some $\tau \in G$
\item
  $ (*): \  \  \Sigma \cap \phi (\Sigma) = \{ 1_G \} .  $
\item
There is no $g \in G$ such that $ \phi(g) \tau g \in \Sigma.$
\end{itemize}

\demo{Acknowledgements}
This paper was written  during visits to  Harvard University and
Florida State University Tallahassee:
I am grateful to both institutions
for their warm hospitality. I would like to thank P. Kronheimer
for asking a good question at the end of my talk in Yau's seminar,
  E. Klassen and V. Kharlamov for a useful conversation, V. Kulikov and
Sandro Manfredini for pointing out some nonsense,  Sandro
again for the nice picture, and finally Grzegorz Gromadzki for pointing
out an error in the second
  version of Lemma~2.3. I want to thank the referee for
helpful comments and Marco Manetti for pointing out that
his examples are not simply connected.
\enddemo

\vglue-8pt
{\it Note.}
Before turning to these examples, I tried to look at rigid surfaces,
trying in particular  to construct nonreal Beauville surfaces.
V. Kharlamov had independently a similar idea, and we spent one
afternoon  together trying to make it work with several examples.
  Later on Kharlamov and Kulikov found the right examples (\cite{K-K}),
one of them a quotient of a  Hirzebruch covering
  of the plane.
 
\vglue-8pt

\AuthorRefNames [Witten]

\end{document}